\documentclass{ieeeconf} 

\usepackage[top=25.4mm,left=19.1mm,right=19.1mm,bottom=19.1mm]{geometry}

\usepackage{graphicx}

\usepackage{xspace}

\usepackage{etoolbox}

\usepackage{amsmath} 
\usepackage{amssymb} 
\usepackage{mathtools}
\usepackage{avldefs}
\usepackage{dsfont} 

\usepackage{algpseudocodex}
\usepackage{algorithm} 

\usepackage{siunitx}
\usepackage{xargs}                      
\usepackage{todonotes}
\usepackage{booktabs}
\usepackage{enumerate}
\usepackage{balance}
\usepackage{tikz}
\usepackage{pgfplots}
\usepackage{gincltex}
\usepackage{pgfmath,pgffor}

\usepackage{xargs}                      
\pgfplotsset{compat=newest} 
\pgfplotsset{
      table/search path={figs/data},
    }
\usetikzlibrary{shapes.geometric,calc,arrows,intersections}

\usepackage[para,online,flushleft]{threeparttable}
\usepackage{multirow}

\makeatletter
\let\NAT@parse\undefined
\makeatother

\usepackage[hyphens]{url}     
\usepackage{hyperref}         
\usepackage{xurl}             
\usepackage{cleveref}

\crefname{ALC@line}{line}{lines}   
\Crefname{ALC@line}{Line}{Lines}   

\usepackage{printlen}

\newcommand{\norm}[1]{\left\lVert#1\right\rVert}
\newcommand{\normm}[2]{{\left\lVert#1\right\rVert}_{#2}}

\DeclareMathOperator*{\argmin}{argmin}

\newtheorem{remark}{Remark}

\makeatletter
\newcommand{\removelatexerror}{\let\@latex@error\@gobble}
\makeatother

\definecolor{mblue}{HTML}{1F77B4}
\definecolor{morange}{HTML}{FF7F0E}
\definecolor{mred}{HTML}{D62728}
\definecolor{mpurple}{HTML}{9467BD}
\definecolor{mbrown}{HTML}{8C564B}
\definecolor{mgreen}{HTML}{00802b}

\definecolor{mosekcolor}{rgb}{0.863,0.129,0.302}
\definecolor{scscolor}{rgb}{0,0.549,0}   
\definecolor{cosmocolor}{rgb}{0,0.392,0.871} 
\definecolor{cosmocdcolor}{rgb}{1,0.4,0} 
\definecolor{osqpcolor}{rgb}{0.4,0.0,0.4} 
\definecolor{gurobicolor}{RGB}{11, 91, 35}

\DeclareMathOperator*{\minimise}{minimise \;\;}

\DeclareMathOperator*{\subjectto}{subject\ to \;\;}

\newcommand{\T}{^\top}

\renewcommand{\Re}{\mathbb{R}}

\let\OLDthebibliography\thebibliography
\renewcommand\thebibliography[1]{
  \OLDthebibliography{#1}
  \setlength{\parskip}{0pt}
  \setlength{\itemsep}{0pt plus 0.3ex}
}

\newcommand{\printColumnWidth}{
    \noindent Column width: \the\columnwidth\ (points)\\[1ex]
    \noindent Column width: \printlength{\columnwidth} (inches)
}

\makeatletter
\newcommand{\linelabel}[1]{%
  \ifcsname c@ALC@line\endcsname
    \protected@edef\@currentlabel{\arabic{ALC@line}}%
  \else
    \protected@edef\@currentlabel{0}%
  \fi
  \label{#1}%
}
\makeatother

\title{\LARGE \bf Krylov Subspace Acceleration for First-Order \\ Splitting Methods in Convex Quadratic Programming}

\author{Gabriel Berk Pereira, Paul J.\ Goulart$^{\dagger}$
\thanks{$^{\dagger}$Department of Engineering Science, University of Oxford. Email: \{\texttt{gabriel.pereira}, \texttt{paul.goulart}\}\texttt{@eng.ox.ac.uk}}
}

\begin{document}
\maketitle


\begin{abstract}
    We propose an acceleration scheme for first-order methods (FOMs) for convex quadratic programs (QPs) that is analogous to Anderson acceleration and the Generalized Minimal Residual algorithm for linear systems.
    We motivate our proposed method from the observation that FOMs applied to QPs typically consist of piecewise-affine operators.
    We describe our Krylov subspace acceleration scheme, contrasting it with existing Anderson acceleration schemes and showing that it largely avoids the latter's well-known ill-conditioning issues in regions of slow convergence.
    We demonstrate the performance of our scheme relative to Anderson acceleration using standard collections of problems from model predictive control and statistical learning applications.
    We show that our method is faster than Anderson acceleration across the board in terms of iteration count, and in many cases in computation time, particularly for optimal control and for problems solved to high accuracy.
\end{abstract}
\section{Introduction} \label{sec:introduction}

We consider throughout the following general convex quadratic program (QP):
\begin{equation} \tag{QP} \label{eq:cone-primal}
    \begin{aligned}
        \minimise &\frac12 x\T P x + c\T x \\
        \subjectto  &A x + s = b \\
        &s \in \Re_+^{m_1} \times \{0\}^{m_2} \\
    \end{aligned}
\end{equation}
with decision variables \(x \in \mathbb{R}^n\) and \(s \in \mathbb{R}^m\), positive semidefinite Hessian matrix \(P \in \mathbb{S}^n_+\), \(c \in \mathbb{R}^n\), \(A \in \mathbb{R}^{m \times n}\), and \(b \in \mathbb{R}^m\).
The feasible set for the slack variable \(s\) encodes \(m_1\) inequality constraints and \(m_2\) equality constraints (\(m_1 + m_2 = m\)).

Problems of the form \eqref{eq:cone-primal} arise in numerous applications in engineering, operations research, finance, signal processing, and other fields \cite{boydConvexOptimization2004a}.
Model predictive control (MPC) applications often require a solution of a QP at each time step \cite{odonoghueSplittingMethodOptimal2013,kouzoupisProperAssessmentQP2015}.
In machine learning, this includes support vector machines \cite{cortesSupportVectorNetworks1995}, Lasso \cite{tibshiraniRegressionShrinkageSelection1996}, and Huber fitting problems \cite{huberRobustEstimationLocation1992}.
Many signal processing problems of interest can also be formulated as QPs \cite{mattingleyRealTimeConvexOptimization2010}.
In finance, they are common in portfolio optimisation \cite{boydMultiPeriodTradingConvex2017}.
Besides these, techniques such as sequential quadratic programming (SQP) for nonlinear programs \cite[§18]{nocedalNumericalOptimization2006a} and branch-and-bound methods for mixed-integer optimisation \cite{fletcherNumericalExperienceLower1998} have QP subproblems at their core.

While interior point methods (IPMs) --- increasingly popular during the 1980s and 1990s --- became standard in many commercial solvers by the turn of the century \cite{wrightPrimalDualInteriorPointMethods1997,potraInteriorpointMethods2000}, since the 2010s there has been a resurgent interest in first-order methods (FOMs) for both very large-scale and embedded optimisation applications \cite{boydDistributedOptimizationStatistical2011,parikhProximalAlgorithms2014a,stellatoOSQPOperatorSplitting2020}.

Broadly speaking, FOMs have much cheaper iterations than IPMs --- which require a matrix factorisation at every step --- but their performance is much more sensitive to problem conditioning \cite{ghadimiOptimalParameterSelection2015,stellatoOSQPOperatorSplitting2020}.
They can be quicker to find modest-accuracy solutions, particularly in large problems for which linear system factorisation can be very expensive, which is tolerable especially in data science applications where the objective function is often a surrogate for another performance metric of interest \cite{boydDistributedOptimizationStatistical2011}.
In embedded applications like MPC, advantages also include the ease with which they can be warm-started and the possibility to make them division-free at deployment \cite{stellatoOSQPOperatorSplitting2020}.

\subsection{FOM Acceleration in Convex Solvers} \label{sec:acceleration-in-solvers}

A large variety of FOM acceleration techniques are used to mitigate potentially slow convergence \cite{daspremontAccelerationMethods2021}.
Nesterov's fast gradient method for smooth unconstrained optimisation \cite[§2]{nesterovLecturesConvexOptimization2018} and the fast iterative shrinkage-thresholding algorithm (FISTA) for non-smooth problems \cite{beckFastIterativeShrinkageThresholding2009} are canonical examples of ``momentum'' acceleration.
However, incorporation of these ideas into general-purpose FOM-based solvers is uncommon.
Analogous schemes have been developed for the widely used alternating direction method of multipliers (ADMM) \cite{boydDistributedOptimizationStatistical2011}, but their success has been limited and restricted to problems satisfying stronger assumptions, such as strong convexity \cite[§1]{zhangGMRESAcceleratedADMMQuadratic2018}.

On the other hand, Anderson acceleration (AA) \cite{andersonIterativeProceduresNonlinear1965} is useful for practical performance improvement in the general-purpose solvers \texttt{SCS} \cite{odonoghueOperatorSplittingHomogeneous2021} and \texttt{COSMO} \cite{garstkaCOSMOConicOperator2021}.
It proposes a candidate iterate by constructing an affine combination of past iterates that minimises the norm of the corresponding combination of past fixed-point residuals \cite[§3]{daspremontAccelerationMethods2021}, and in large conic programs (including QPs) can be used at little additional cost.
However, it suffers from conditioning problems when past fixed-point residuals are nearly collinear, which frequently happens in regions of slow convergence \cite{garstkaSafeguardedAndersonAcceleration2022}.   The solver
\texttt{SCS} seeks to avoid this by only updating its AA module every 10 iterations by default \cite{scs}.
\texttt{COSMO} rejects the candidate outright if the norm of the inner least squares solution is too large, indicating a poorly-conditioned subproblem \cite{garstkaSafeguardedAndersonAcceleration2022}.

Using the fact that ADMM for equality-constrained QPs consists of linear updates, the authors in \cite{zhangGMRESAcceleratedADMMQuadratic2018} apply the Generalized Minimal Residual (GMRES) method \cite{saadGMRESGeneralizedMinimal1986}, originally devised to solve linear systems, to this context.
They briefly observe that whenever iterate update equations are nonlinear (\eg ADMM for general QPs), they ``might still be approximated by their linearisation'', but do not further pursue this direction.

In fact, there is a close equivalence between AA and GMRES in iterative methods with affine dynamics \cite[Thm.~2.2]{walkerAndersonAccelerationFixedPoint2011}.
The latter includes, as a relevant example, ADMM applied to \textit{equality-constrained} QPs, an observation which authors exploit in \cite{zhangGMRESAcceleratedADMMQuadratic2018}.
Our proposed acceleration method leverages a local linearisation of the FOM operator, which for QPs is exact within active-set regions due to the piecewise-affine structure of the dynamics.
The resulting acceleration procedure is nearly equivalent to Anderson acceleration in exact arithmetic (in the sense of \cite{walkerAndersonAccelerationFixedPoint2011}) but sidesteps its conditioning problems through the Arnoldi process's orthogonalisation \cite[§33]{trefethenNumericalLinearAlgebra1997}.

\subsection{Paper Outline}

In \Cref{sec:foms-on-qps}, we recall notions relating to FOMs in convex optimisation as well as the well-known ADMM baseline which we use thereafter for illustration.
In \Cref{sec:krylov-acceleration}, we motivate and propose a novel Krylov subspace acceleration procedure which is equivalent in a sense to Anderson acceleration, while largely sidestepping the latter's ill-conditioning issues.
Then, in \Cref{sec:numerical-experiments}, we compare the numerical performance of our Krylov method against Anderson acceleration on 104 benchmark problems arising in model predictive control and statistical learning.

All reported numerical experiments were performed single-threaded using \texttt{Julia v1.11.5} on an Intel Xeon w9-3475X workstation with 256GB RAM.

\subsection{Notation}

We denote the identity matrix of dimension \(d\) with \(I_d\), or \(I\) where the dimension is clear.
Subscripts denote iteration counters.
For vectors, they may also denote entry indices; context clarifies the use case.
The symbol \(\mathbb{S}^d_+\) denotes the set of symmetric positive semidefinite matrices of dimension \(d \times d\).
We write \(A \succeq 0\) when \(A\) is positive semidefinite.

The indicator function of a set \(C\) is \(\mathcal{I}_{C}\), valued \(\infty\) for inputs in \(C\) and \(0\) otherwise.
The projection operator onto a set \(C\) with respect to the Euclidean norm is \(\Pi_C\).
The norm with which a normed vector space is endowed is denoted by \(\norm{\cdot}\).
The common Euclidean norm is \(\normm{\cdot}{2}\).
\section{First-Order Methods for Quadratic Programming} \label{sec:foms-on-qps}

FOMs for convex optimisation are often posed as iterations of some  fixed-point operator denoted by \(T\), for which two properties are sought.
First, that the set of its fixed points is the set of solutions to the optimisation problem of interest.
Second, that it has properties, \eg averagedness or nonexpansiveness, which allow us to straightforwardly design iterative algorithms whose iterates converge to that fixed point set.
In simple cases, an iteration of the form \(u_{k+1} \gets T u_k\) converges to a solution to the optimisation problem \cite{bauschkeConvexAnalysisMonotone2017,ryuLargeScaleConvexOptimization2022}.
We reserve the symbol \(T\) to generically denote the operator defining a FOM under consideration.

This paper's contribution relies on the observation that FOMs for QPs often consist of a piecewise-affine (PWA) iterate update operator.
We are not the first to observe this property, and we refer the reader to \cite[§2]{ranjanVerificationFirstorderMethods2025} for a description in these terms of a suite of primitive algorithmic steps found in this setting.

We will use the well-known alternating direction method of multipliers (ADMM) \cite{boydDistributedOptimizationStatistical2011} as an illustrative example here, as well as in the benchmarks of \Cref{sec:numerical-experiments}, though our approach is generically applicable to any FOM that amounts to a PWA mapping.
In various guises, ADMM is the basis for such convex solvers as \texttt{OSQP} \cite{stellatoOSQPOperatorSplitting2020}, \texttt{SCS} \cite{odonoghueOperatorSplittingHomogeneous2021}, and \texttt{COSMO} \cite{garstkaCOSMOConicOperator2021}.
For this, we recast the problem \eqref{eq:cone-primal} into one with a composite objective and linear constraints.
Let \(f(x) \coloneq \frac12 x\T P x + c\T x\) and \(g(s) \coloneq \mathcal{I}_{\mathcal{K}}(s)\), where \(\mathcal{K} \coloneq \Re_+^{m_1} \times \{0\}^{m_2}\).
Then we drop the inclusion constraints in \eqref{eq:cone-primal} and let the objective of the recast problem be \(f(x) + g(s)\).
A common formulation of ADMM \cite[§3]{boydDistributedOptimizationStatistical2011} is then
\begin{align}
    x_{k+1} &= \argmin_{x \in \Re^n} \left\{f(x) - \langle y_k, A x \rangle + \frac\rho2 \normm{b - A x - s_k}{2}^2\right\} \notag \\
    s_{k+1} &= \argmin_{s \in \Re^m} \left\{ g(s) - \langle y_k, s \rangle + \frac\rho2 \normm{b - Ax_{k+1} - s}{2}^2 \right\} \notag \\
    y_{k+1} &= y_k + \rho \left(b - A x_{k+1} - s_{k+1} \right). \notag
\end{align}

It is convenient for us to manipulate this a bit before proceeding.
We start by eliminating the slack variable \(s\) from the updates above with elementary algebra and simplifying expressions by substituting in \(f\) and \(g\).
If (and only if) \(P + \rho A\T A \nsucc 0\), we also add a regularisation term \(\frac{\delta}{2} \normm{x - x_k}{2}^2\) to the \(x\) minimisation subproblem objective, with small \(\delta > 0\).\footnote{This ensures that the minimiser \(x_{k+1}\) exists and is unique. Asymptotic convergence of the method is retained for any \(\delta > 0\); see \cite[§3.2]{shefiRateConvergenceAnalysis2014}.}
In general we let \(W \coloneq P + \rho A\T A + \delta I\), with nonnegative \(\delta\).
The resulting iteration is
\begin{subequations} \label{eq:qp-preppm-updates}
    \begin{align}
        y_{k+1} &= \Pi_{\mathcal{K}^*} \left(y_k + \rho (A x_k - b)\right) \label{eq:qp-preppm-y-update} \\
        \bar{y}_{k+1} &= 2 y_{k+1} - y_k \\
        x_{k+1} &= x_k - W^{-1} \left(P x_k + c + A^T \bar{y}_{k+1} \right). \label{eq:qp-preppm-x-update}
    \end{align}
\end{subequations}
The dual cone in \eqref{eq:qp-preppm-y-update} is \(\mathcal{K}^* \coloneq \Re^{m_1}_+ \times \Re^{m_2}\).
That is, the multipliers of the \(m_1\) inequality constraints are projected to the nonnegative cone, while those of the \(m_2\) equality constraints are ``free''.
Note that the term \(\bar{y}_k\) is only for simplicity of exposition --- the iteration is carried by \((x, y)_{k+1} \gets T (x, y)_k \in \Re^{m + n}\).

It is straightforward to verify that the preceding is a preconditioned proximal point method \cite{brediesDegeneratePreconditionedProximal2022} applied to a standard saddle-point formulation of \eqref{eq:cone-primal}, with preconditioning matrix
\begin{equation} \label{eq:admm-preconditioner}
    M = 
    \begin{bmatrix}
        \rho A\T A + \delta I & A\T \\
        A & \frac1\rho I_m
    \end{bmatrix} \succ 0;
\end{equation}
see \cite[§3.5]{chambolleIntroductionContinuousOptimization2016} and \cite[§3.2]{liuAccelerationPrimalDual2021}.
The proximal point method \cite{rockafellarMonotoneOperatorsProximal1976} is commonly used in analysis of algorithms of this type, so this form will be useful when designing a safeguard for our acceleration method.

\section{Krylov Acceleration of FOMs for QPs} \label{sec:krylov-acceleration}

\subsection{Local Linearisation of FOM Operators} \label{sec:fom-linearisation}

Given a FOM operator \(T\) as in \Cref{sec:foms-on-qps}, linearising the fixed-point equation \(T(u) = u\) around a current iterate \(u_k \in \Re^d\) yields
\begin{equation} \label{eq:linearised-fp-system}
    (L - I) \delta = -r_k,
\end{equation}
where \(L \coloneq DT(u_k)\) is the Jacobian of \(T\) at \(u_k\), \(r_k \coloneq T u_k - u_k\) is the fixed-point residual, and \(\delta\) is a correction to \(u_k\) approximating a fixed point of \(T\).

For QPs, the operator \(T\) is piecewise-affine (PWA), so \(L\) is piecewise-constant and the linearisation \eqref{eq:linearised-fp-system} is exact within each affine region.
To see this, consider the updates for \((x, y)_k\) in \eqref{eq:qp-preppm-updates}.
The primal update \eqref{eq:qp-preppm-x-update} is an affine function of \((x, y)_k\).
The dual update \eqref{eq:qp-preppm-y-update} is a composition of an affine function and the projection \(\Pi_{\mathcal{K}^*}\).
While \((x, y)_{k+1}\) is therefore generally a non-affine function of \((x, y)_k\), this projection operator is still simple in some sense.
Indeed, \(\Re^{m + n}\) can be partitioned into polyhedra within each of which the transformation describing the primal-dual iterate update is locally affine.

Consider the following example.
We let \(a_j\T\) denote the \(j\)th row of \(A\).
Given some \(x_k \in \Re^n\), within the polyhedron in \(\Re^{m}\) defined by the system of \(m_1\) inequalities
\begin{equation}
    \begin{aligned}
        y_1 + \rho \left( a_1\T x_k - b_1 \right) &\leq 0 \\
        y_2 + \rho \left( a_2\T x_k - b_2 \right) &\geq 0 \\
        \cdots \\
        y_{m_1} + \rho \left( a_{m_1}\T x_k - b_{m_1} \right) &\geq 0,
    \end{aligned}
\end{equation}
the projection operation in \eqref{eq:qp-preppm-y-update} is linear: it acts as the identity for entries 2 to \(m\), and as the constant zero function for the first entry.
Such a system of inequalities defines an active set \(\mathcal{J} \subseteq \{1, \ldots, m_1\}\).
Here, \(\leq 0\) means that the corresponding pre-projection argument is nonpositive, so the projection zeros it out (inactive constraint, zero multiplier), while \(\geq 0\) means it is nonnegative and the projection acts as the identity (active constraint, nonzero multiplier).
It becomes clear that \(\Re^{m + n}\) can be partitioned into \(2^{m_1}\) polyhedra, within each of which the mapping \eqref{eq:qp-preppm-updates} from \((x, y)_k\) to \((x, y)_{k+1}\) is affine.

Since other FOMs for QPs share a PWA structure \cite[§2]{ranjanVerificationFirstorderMethods2025}, we refer to a generic optimisation iterate \(u \in \Re^d\) for the remainder of this section.
Within each active-set region \(\mathcal{J}\), the Jacobian \(L\) is constant and the linearisation \eqref{eq:linearised-fp-system} is exact.
In our ADMM example, for a given point \(u_k = (x, y)_k\), the active set \(\mathcal{J}\) is determined by the signs of the entries of \(y_k + \rho \left(A x_k - b\right) \).
In other FOMs, the details of this relation may differ.

Note that if the current active set \(\mathcal{J}\) happens to be that of a minimiser --- an ``optimal active set'' --- then the linearisation is exact and \(u_k + \delta^*\) from \eqref{eq:linearised-fp-system} is a QP solution.

\subsection{Krylov Subspace FOM Acceleration}

In \cite{zhangGMRESAcceleratedADMMQuadratic2018}, the authors apply ADMM to a class of equality-constrained QPs.
In this scenario the active set is fixed, so solving \eqref{eq:linearised-fp-system} returns solutions to the optimisation problem.
Leveraging this fact, they use the affine mapping defining their baseline method in a GMRES routine \cite{saadGMRESGeneralizedMinimal1986}, with residual-minimising iterates drawn from nested Krylov subspaces of growing dimension.

In fact, GMRES is intimately related to the original (``type-II'') AA (see \cite[§2]{zhangGloballyConvergentTypeI2020} for an explanation of this as well as ``type-I'' AA).
The authors in \cite{walkerAndersonAccelerationFixedPoint2011} show that, under mild assumptions, there is an immediate connection between AA applied to finding the fixed point of an affine mapping and the GMRES iterates for solution of a closely, naturally related linear system \cite{saadGMRESGeneralizedMinimal1986}.

We use our earlier observations about the ADMM mapping \eqref{eq:qp-preppm-updates} being PWA to extend the method of \cite{zhangGMRESAcceleratedADMMQuadratic2018} to the much broader class of general QPs.
This is trivial if one knows \textit{a priori} an optimal active set of the QP at hand, so that we can obtain solutions to it from an easily reformulated equality-constrained problem.

Generally, of course, FOMs can expend a considerable number of iterations ``searching for'' an optimal active set in the first place \cite{boleyLocalLinearConvergence2013}.
In any case, we consider applying a GMRES routine to the linearised system \eqref{eq:linearised-fp-system} using the current Jacobian \(L\), even though it will generally not correspond to an optimal active set.
The rationale is that we replace the nonlinear iterate update mapping with a local linearisation of the FOM operator, which for QPs is exact within each active-set region.

We perform the Arnoldi process of GMRES \cite[§35]{trefethenNumericalLinearAlgebra1997} in parallel with conventional fixed-point operator updates to the true optimisation iterate.
So long as the Jacobian \(L\) remains invariant for a number of consecutive iterations (which for QPs means an invariant active set \(\mathcal{J}\)), this is a bona fide Arnoldi process for the respective linear operator.
Even if \(L\) changes, we blindly carry out the same steps, resulting in a perturbed Arnoldi recurrence.
A full analysis of the Arnoldi process under such perturbations is beyond the scope of this work; we rely on a safeguarding measure described in \Cref{sec:krylov-adapted-safeguard} to ensure that perturbed candidates do not compromise convergence.

We approximate the solution to \eqref{eq:linearised-fp-system} using GMRES \cite{saadGMRESGeneralizedMinimal1986}, which finds the best approximation (in the sense of minimising the residual's Euclidean norm) within a Krylov subspace.
The method computes an orthogonal basis of a Krylov subspace using the Arnoldi process, which requires computing products through \(L\) as well as the orthogonalisation of basis vectors using the modified Gram-Schmidt (MGS) process \cite[§8]{trefethenNumericalLinearAlgebra1997}.
This circumvents ill-conditioning problems which Anderson acceleration suffers from when recent residuals are nearly collinear \cite[§3B]{garstkaSafeguardedAndersonAcceleration2022}, as noted in \Cref{sec:acceleration-in-solvers}.

\begin{remark} \label{rem:krylov-shift-invariance}
    The system \eqref{eq:linearised-fp-system} has matrix \(L - I\), but the Krylov subspace \(\mathcal{K}_j(L - I, r_k) = \mathcal{K}_j(L, r_k)\) for all \(j\), since Krylov subspaces are invariant under identity shifts of the generating matrix (by the binomial theorem).
    We may therefore equivalently build the Arnoldi basis using \(L\) as the operator, which we do hereafter.
\end{remark}

As in restarted GMRES \cite[§3.1]{saadGMRESGeneralizedMinimal1986}, we cap the dimension of the Krylov subspaces used at some ``memory size'' \(\mu\), then restart (\ie discard Krylov subspace data from previous iterates and start over).
As usual in GMRES, we take the initial basis vector each time to be the system residual given our best estimate of the solution, namely the current optimisation iterate.
At step \(j\) of the Arnoldi process, assuming that the Jacobian has remained invariant since the last restart, this translates to the relation
\begin{equation} \label{eq:arnoldi-relation}
    L Q_{j-1} = Q_j \tilde{H}_{j - 1},
\end{equation}
where \(Q_j \in \Re^{d \times j}\) is a matrix of \(j\) orthonormal Krylov basis vectors (\(Q_{j-1}\) is formed from it by removing the final column) and \(\tilde{H}_{j - 1} \in \Re^{j \times (j - 1)}\) is upper Hessenberg.
But note that we follow these steps even if the Jacobian \textit{does} change, which leads to a (possibly heavily) perturbed analogue of \eqref{eq:arnoldi-relation}.

While this Arnoldi process grows a basis and a Hessenberg matrix at every solver iteration, we can avoid incurring the costs of computing a candidate point at every step.
If we choose to do it at iteration \(k\), with Arnoldi process at step \(j\), we start by solving
\begin{align}
    z_{\text{Kr}} = &\argmin_{z \in \Re^{j-1}} \normm{(L - I) Q_{j-1} z + r_k}{2}^2 \notag \\
    = &\argmin_{z \in \Re^{j-1}} \normm{\left(\tilde{H}_{j-1} - \begin{bmatrix} I_{j-1} \\ 0 \end{bmatrix}\right) z + Q_j\T r_k}{2}^2, \label{eq:hessenberg-least-squares}
\end{align}
where the second equality uses \Cref{rem:krylov-shift-invariance} and the Arnoldi relation \eqref{eq:arnoldi-relation}.
This is a \(j \times (j - 1)\) least squares problem whose Hessenberg structure makes it cheaply solvable through a QR method based on first triangularising \(\tilde{H}_{j-1}\) with \(j - 1\) Givens rotations \cite[§3.2]{saadGMRESGeneralizedMinimal1986}, which we compute and store at each Arnoldi iteration.

From this subproblem we obtain \(u_{\text{Kr}} = u_k + Q_{j-1} z_{\text{Kr}}\).
To reflect the close relationship between type-II AA and GMRES in the sense of \cite[Thm.~2.2]{walkerAndersonAccelerationFixedPoint2011}, we define the Krylov-accelerated candidate as \(\hat{u}_k \coloneq T u_{\text{Kr}}\).

\subsection{Safeguarding} \label{sec:krylov-adapted-safeguard}

The conic solvers \texttt{SCS} and \texttt{COSMO} use safeguards for their AA modules based on the (Euclidean) norm of the fixed-point residual \cite{scs,garstkaCOSMOConicOperator2021}.
Thus they reject candidate accelerated iterates that would  increase the fixed-point residual norm.
Here, we adapt and simplify the safeguarding framework designed for learnt optimisers in \cite{heatonSafeguardedLearnedConvex2023}.
The authors there consider only FOMs which are averaged in the Euclidean sense, but one can straightforwardly extend the analysis to averagedness in general metrics.

Suppose we have a FOM operator \(T\) which is averaged in a norm \(\norm{\cdot}\).
At arbitrary iterations \(k\), an arbitrary acceleration module proposes a point \(\hat{u}_k\).
Say we impose the safeguard
\begin{align}
    \norm{T \hat{u}_k - \hat{u}_k} &\leq \eta \norm{T u_k - u_k} \text{ and} \label{eq:residual-contraction} \\
    \norm{\hat{u}_k - u_k} &\leq c \norm{T u_k - u_k} \label{eq:step-size-bound}
\end{align}
with parameters \(\eta \in (0, 1)\) and \(c > 0\) as a condition for assigning \(u_{k+1} \gets T \hat{u}_k\).
If rejected, a standard (non-accelerated) iteration \(u_{k+1} \gets T u_k\) is taken instead.
Then the accelerated method retains its basic convergence properties in the sense of \cite[Thm.~2]{heatonSafeguardedLearnedConvex2023}.

Preconditioned proximal point methods are averaged in the norm induced by the preconditioning matrix \cite[Lem.~2.6]{brediesDegeneratePreconditionedProximal2022}.
Thus, for the ADMM updates \eqref{eq:qp-preppm-updates}, the appropriate norm in which to evaluate the safeguard is \(\normm{\cdot}{M}\) with \(M\) as in \eqref{eq:admm-preconditioner}.

A summary of the safeguarded Krylov-accelerated method we have described is presented in \Cref{alg:krylov-acc-optimisation}.
In this pseudocode, we reserve matrix subscripts to index working matrices.
Namely, given a matrix \(A\), \(A_{i, j}\) denotes the entry in row \(i\), column \(j\).
We use colons (:) to denote index ranges defining submatrices and subvectors in the natural way.

\begin{algorithm}[ht]
    \caption{Safeguarded Krylov-accelerated FOM for QP}
    \label{alg:krylov-acc-optimisation}
    \begin{algorithmic}[1]
        \Require Baseline FOM \(\alpha\)-averaged operator \(T\) with Jacobian \(L = DT(u)\). Safeguard parameters \(\eta \in (0, 1)\) and \(c > 0\). Memory size \(\mu\). Acceleration attempt set \(\mathcal{G} \subseteq \{3, \ldots, \mu + 1\}\) with \(\mu + 1 \in \mathcal{G}\).
        \State \(k \gets 0, j \gets 1, H \gets 0^{(\mu + 1) \times \mu}, Q \gets 0^{d \times \mu}\)
        \While{not converged}
            \If{\(j \in \mathcal{G}\)} \Comment{Attempt accelerated step}
                \State Solve \eqref{eq:hessenberg-least-squares} for \(z_{\text{Kr}} \in \Re^{j-1}\)
                \State \(u_{\text{Kr}} \gets u_k + Q_{1:d, 1:j-1} z_{\text{Kr}}\)
                \State \(\hat{u}_k \gets T u_{\text{Kr}}\) \label{line:krylov-operator-step}
                \If{\(\hat{u}_k\) satisfies \eqref{eq:residual-contraction} and \eqref{eq:step-size-bound}} \label{line:safeguard-check}
                    \State \(u_{k+1} \gets T \hat{u}_k\)
                \EndIf
                \If{\(j > \mu\)} \Comment{Memory full, restart}
                    \State \(H \gets 0^{(\mu + 1) \times \mu}, Q \gets 0^{d \times \mu}, j \gets 0\) \label{line:gmres-restart}
                \EndIf
            \Else
                \State \(u_{k+1} \gets T u_k\)
                \State Compute \(L = DT(u_k)\)
                \If{\(j = 1\)} \Comment{Initialise Krylov basis}
                    \State \(q \gets r_k / \normm{r_k}{2}\), where \(r_k \coloneq T u_k - u_k\)
                    \State \(Q_{1:d, 1} \gets q\)
                \Else \Comment{Arnoldi step}
                    \State \(q \gets L q\)
                    \State Orthogonalise \(q\) against \(Q_{1:d, 1:j-1}\) with MGS and normalise, in-place. Store coefficients in \(H_{1:j, j-1}\).
                    \State \(Q_{1:d, j} \gets q\)
                    \State Apply and store Givens rotation to eliminate new subdiagonal entry \(H_{j, j-1}\) in-place.
                \EndIf
            \EndIf
            \State \(k \gets k + 1, j \gets j + 1\)
        \EndWhile
    \end{algorithmic}
\end{algorithm}

\subsection{Computational Considerations}

\subsubsection{Arnoldi Process Concurrency}

The reader might suppose that an iteration of \Cref{alg:krylov-acc-optimisation} is much more expensive than a standard FOM iteration, since the Jacobian-vector product \(Lq\) required by the Arnoldi step at each iteration seemingly incurs practically the same costs as an application of the FOM operator \(T\) itself.

But in marginal terms this is not so, since \(L\) shares most of its computation with the operator application \(Tu\).
Consider keeping a complex working vector\footnote{or, more generally, any vector composed of two-element data structures.} whose real part is the optimisation iterate \(u_k\) and imaginary part is the Arnoldi working vector denoted by \(q\) in \Cref{alg:krylov-acc-optimisation}.
The costliest primitives in the application of a typical FOM operator for QPs, such as \eqref{eq:qp-preppm-updates} --- vector products through sparse matrices; implicit multiplication by a matrix inverse; vector addition --- can be performed on this complex working vector, with significant cuts to the marginal cost of these typically memory-bound operations.

The exception is the projection, in our case to \(\mathcal{K}^* = \Re^{m_1}_+ \times \Re^{m_2}\).
We must first perform it for the optimisation iterate, deducing the Jacobian in the process.
Then we can apply its linearisation to the Arnoldi working vector, and proceed thereafter processing both simultaneously again.
Since this is a very cheap primitive, the serial duplication of work does not constitute a great performance loss.

\subsubsection{Primitive Operation Costs} \label{sec:primitive-costs}

Accounting for the fact that many QP application problems --- \eg in statistical learning, control, and portfolio optimisation --- have sparse data \cite{stellatoOSQPOperatorSplitting2020}, matrices \(P\) and \(A\) are stored and manipulated in our code in Compressed-Sparse-Column (CSC) format \cite{davisDirectMethodsSparse2006} (as in \texttt{OSQP}, for instance).
We are therefore concerned with computational costs in the context of CSC matrix representations.

One crucial such matter is that of sparse matrix-vector multiply (SpMV) operations \cite{williamsOptimizationSparseMatrixVector2007}.
These have low operational intensity, \ie few flops are executed per byte of retrieved memory.
Hence the limiting factor in SpMV speed tends to be memory, not floating-point performance, \ie this is a memory-bound operation \cite{williamsRooflineInsightfulVisual2009}.
Similar considerations apply when solving linear systems given sparse factor structures (in our case, a sparse Cholesky factorisation \cite[§3]{davisDirectMethodsSparse2006}).

Thus, one can exploit memory locality to perform two of these operations simultaneously in less than twice the time it takes to compute one.
As we have mentioned already, a simple way to achieve this is to compute these products with complex vectors.

Using \texttt{@belapsed} from the \texttt{BenchmarkTools.jl} library in \texttt{Julia}, we benchmarked SpMV and in-place Cholesky linear solutions --- with coefficient matrices \(W\) as in \eqref{eq:qp-preppm-updates} --- using vectors with (64-bit) real as well as complex entries.
We used \(P\) and \(A\) matrices from subsets of the \texttt{mpc} and \texttt{sslsq} problem collections which we describe in \Cref{sec:numerical-experiments}.
For each of the SpMV operations through \(P\), \(A\), and \(A\T\), with a real and a complex vector, we take the minimum time in the \texttt{@belapsed} call as representative.

\begin{figure}[ht]
    \centering
    \includegraphics{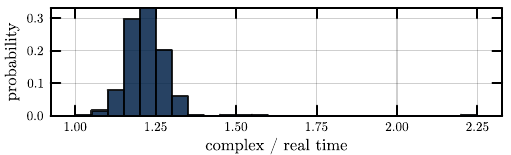}
    \caption{Histogram of ratio between SpMV (\(P x\), \(A x\), and \(A\T y\)) times with complex/real 64-bit floating point vectors. CSC matrices \(P\) and \(A\) taken from 108 \texttt{mpc} and \texttt{sslsq} problems with \(m, n \leq \num{100000}\). Ratios are concentrated well below 2, confirming that the marginal cost of a simultaneous second SpMV is small. (3 samples with ratio below 1, attributed to measurement variability, were removed.)}
    \label{fig:spmv-real-complex-hist}
\end{figure}

\begin{figure}[ht]
    \centering
    \includegraphics{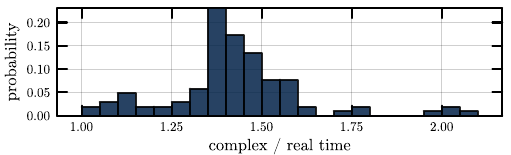}
    \caption{Histogram of ratio between time to solve linear systems in-place with matrix \(W\) for one versus two right-hand sides. \(W\) as in \eqref{eq:qp-preppm-updates} with \(\rho = 0.1\). CSC matrices \(P\) and \(A\) taken from 104 \texttt{mpc} and \texttt{sslsq} problems with \(m, n \leq \num{20000}\). As with SpMV, ratios are well below 2.}
    \label{fig:chol-solve-real-complex-hist}
\end{figure}

The complex-to-real time ratios of these operations on these data are shown in histograms in \Cref{fig:spmv-real-complex-hist,fig:chol-solve-real-complex-hist}.

\section{Numerical Benchmarks} \label{sec:numerical-experiments}

In this section we present the results of numerical experiments with the purpose of comparing safeguarded (type-II) Anderson- and Krylov-based acceleration methods for FOMs on QPs.
The benchmark problem sets we consider are among those used in \cite[§5]{goulartClarabelInteriorpointSolver2024}, which we refer the reader to for more details on QP formulations and applications.

\subsection{Benchmark Setup}

For a baseline method we use ADMM as in \eqref{eq:qp-preppm-updates}.
We fix \(\rho = 0.1\) throughout (the starting step size in \texttt{OSQP} \cite[§5.2]{stellatoOSQPOperatorSplitting2020}).
As mentioned in the discussion of \eqref{eq:qp-preppm-updates}, wherever the Cholesky routine provided by \texttt{SparseArrays.jl} fails at first with \(\delta = 0\), we set \(\delta = 10^{-10}\) for positive definiteness of \(W\).
For both Anderson- and Krylov-accelerated methods, we use memory size \(\mu = 15\).

When evaluating the safeguard, we only check condition \eqref{eq:residual-contraction}, omitting the step-size bound \eqref{eq:step-size-bound}.
Since rejected candidates default to a standard FOM step, using any \(\eta < 1\) ensures that the residual norm sequence converges to zero asymptotically \cite[Thm.\ 1]{giselssonLineSearchAveraged2016}.
For simplicity, we use \(\eta = 1\); this is consistent with the safeguard used in \texttt{SCS} \cite{scs}.

We use two Krylov variants differing on the frequency with which the method proposes an accelerated iterate, \ie differing \(\mathcal{G}\) in \Cref{alg:krylov-acc-optimisation}.
In one, we have \(\mathcal{G} = \{16\}\).
This means that we only carry out the steps to propose an accelerated point once the Krylov basis matrix has been fully populated with 15 orthogonal vectors.
In the other, we use \(\mathcal{G} = \{6, 11, 16\}\), proposing a new accelerated iterate 3 times per restart.
The trade-off is as follows.
Proposing accelerated points more frequently can be costlier, owing especially to additional applications of \(T\) in evaluating the candidate as well as the safeguard. 
However, this can mitigate issues caused by active set changes during the Arnoldi process, which perturb relation \eqref{eq:arnoldi-relation} and often render subsequent candidates before the next restart unusable.\footnote{This is not always the case: we occasionally observe accelerated points of good quality following a perturbed Arnoldi process. Perturbation analysis of Krylov methods can be found in \cite{simonciniTheoryInexactKrylov2003}.}
In figures we denote these variants with ``\(\text{tries}=1\)'' and ``\(\text{tries}=3\)'' respectively.\footnote{Memory size 15 is large enough to explore the use of several Krylov acceleration attempts per restart, without becoming too vulnerable to stale Jacobians. For reference, \texttt{SCS} and \texttt{COSMO} use 10 by default.}

We use an efficient type-II AA implementation with restarted memory and QR-based least squares solution, provided in the \href{https://github.com/oxfordcontrol/COSMOAccelerators.jl}{\texttt{COSMOAccelerators.jl}} package and used in the \texttt{COSMO} solver.
This can mostly be thought of as plugged into \Cref{alg:krylov-acc-optimisation} in lieu of our proposed Krylov acceleration steps, but a couple of differences are worth noting.
For one, we let the type-II AA method propose an iterate every time the accelerator is updated with a new input/output pair (unlike in our Krylov method, this alone does \textit{not} require applying \(T\)).
Besides the obvious analogue, we also use a variant of AA which only lets the acceleration module ``see'' every 10\textsuperscript{th} iterate.
This is equivalent to applying AA to the composition \(T^{10}\), and it is used by default in \texttt{SCS} as a way to mitigate numerical issues associated with near collinearity of residuals \cite{scs}.
In figures we denote these variants with ``\(\text{interval}= 1\)'' and ``\(\text{interval}= 10\)'' respectively.

We run each combination of solver and problem twice and take the best wall-time performance to represent it.

\subsection{Performance Profiles}

We monitor solution progress using the following residual measures, which we update every 25 iterations:
\begin{align}
    r_p &= \frac{\normm{[A x - b]_+}{\infty}}{1 + \max\{\normm{A x}{\infty}, \normm{b}{\infty}\}} \notag \\
    r_d &= \frac{\normm{P x + A\T y + c}{\infty}}{1 + \max \{\normm{P x}{\infty}, \normm{A\T y}{\infty}, \normm{c}{\infty}\}} \notag \\
    \text{pd} &= \frac{\lvert x\T P x + c\T x + b\T y\rvert}{1 + \max \{\lvert (1/2) x\T P x + c\T x \rvert, \lvert (1/2) x\T P x + b\T y \rvert\}} \notag.
\end{align}
Here, \([\cdot]_+\) denotes a vector's entry-wise nonnegative part.
We declare a problem solved with accuracy \(\epsilon > 0\) when \(\max \left\{r_p, r_d, \text{pd}\right\} \leq \epsilon\) (this is also done \eg in \cite[§5]{luPracticalOptimalFirstOrder2025}).

We compare methods using (relative) performance profiles \cite{dolanBenchmarkingOptimizationSoftware2002}.
Given a set of test problems, for each solver \(s\) and problem \(p\) we let \(t_{p, s}^{\epsilon}\) denote the budget (either iteration count or wall time) used by solver \(s\) in solving problem \(p\) to accuracy \(\epsilon\) (set to \(\infty\) if the solver fails to solve the problem within \(\num{20000}\) iterations).
We then define the relative performance ratio as
\begin{equation}
    u_{p, s}^{\epsilon} = \frac{t_{p, s}^{\epsilon}}{\min_s t_{p, s}^{\epsilon}}.
\end{equation}

In a given benchmark, we let \(N_{\text{solved}}\) denote the number of problems which were solved by at least one of the solvers under consideration.\footnote{This differs from the normalisation originally used in \cite{dolanBenchmarkingOptimizationSoftware2002}, namely the total number of test problems.}
The relative performance profile for a solver \(s\) at accuracy \(\epsilon\) plots the function \(P^{\epsilon}_s \colon [1, \infty) \mapsto [0, 1]\) which we define as
\begin{equation}
    P^{\epsilon}_s(\tau) = \frac{1}{N_{\text{solved}}} \sum_p \mathds{1}_{\{\cdot \leq \tau\}} (u_{p, s}^{\epsilon}),
\end{equation}
where \(\mathds{1}_{\{\cdot \leq \tau\}} (u)\) equals 1 if \(u \leq \tau\) and 0 otherwise.
All profiles use a logarithmic horizontal axis.

\subsection{Model Predictive Control Problems \textup{(\texttt{mpc})}}

\begin{figure}[ht]
    \centering
    \includegraphics{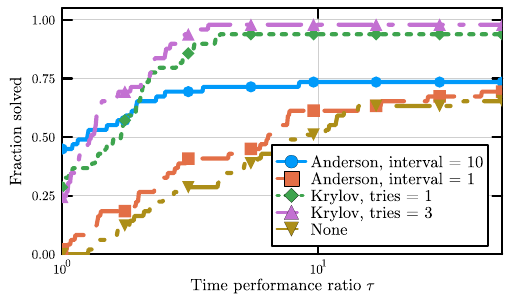}
    \caption{Relative time performance profile for tolerance \(\epsilon = 10^{-3}\) on subset of \texttt{mpc} problems used in \cite{goulartClarabelInteriorpointSolver2024}. 49 out of 72 problems were solved by at least one solver.}
    \label{fig:admm-time-relative-mpc-eps1e-3}
\end{figure}

\begin{figure}[ht]
    \centering
    \includegraphics{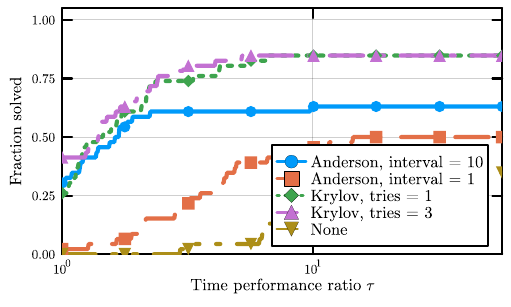}
    \caption{Relative time performance profile for tolerance \(\epsilon = 10^{-6}\) on subset of \texttt{mpc} problems used in \cite{goulartClarabelInteriorpointSolver2024}. 46 out of 72 problems were solved by at least one solver.}
    \label{fig:admm-time-relative-mpc-eps1e-6}
\end{figure}

We use 72 QPs derived from MPC problems in industrial and academic settings from \cite{kouzoupisProperAssessmentQP2015}, with \(m, n \leq \num{20000}\).
These are the same MPC instances as used in the paper presenting the \texttt{Clarabel} IPM solver \cite{goulartClarabelInteriorpointSolver2024}.
Time profiles for \(\epsilon = 10^{-3}\) and \(\epsilon = 10^{-6}\) are shown in \Cref{fig:admm-time-relative-mpc-eps1e-3,fig:admm-time-relative-mpc-eps1e-6} respectively.

The Krylov methods broadly outperform the Anderson ones on this benchmark set, with 3 tries per memory yielding a wall time advantage over the other Krylov variant.
The Anderson variant applied to \(T^{10}\) does better than that which is applied to \(T\) itself, demonstrating the potential benefit of this technique in mitigating conditioning issues in least squares subproblems.

Iteration-count profiles (omitted for space) show a consistent advantage for the Krylov variants across both tolerance levels.



\subsection{SuiteSparse Least Squares Problems \textup{(\texttt{sslsq})}}

\begin{figure}[ht]
    \centering
    \includegraphics{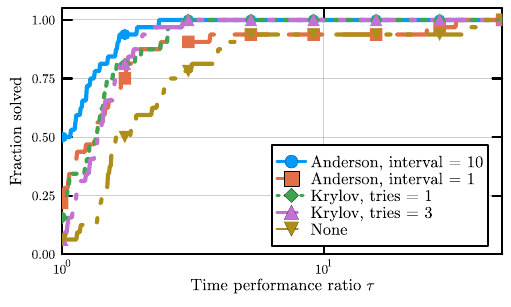}
    \caption{Relative time performance profile for tolerance \(\epsilon = 10^{-3}\) on subset of \texttt{sslsq} problems with \(m, n \leq \num{20 000}\). All 32 problems were solved by at least one solver.}
    \label{fig:admm-time-relative-sslsq-eps1e-3}
\end{figure}

\begin{figure}[ht]
    \centering
    \includegraphics{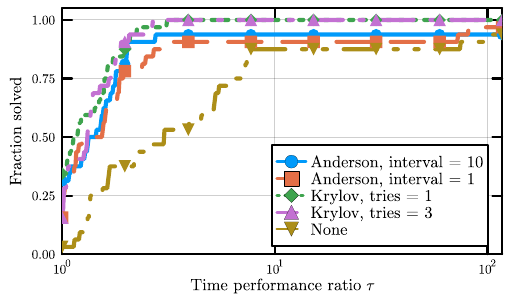}
    \caption{Relative time performance profile for tolerance \(\epsilon = 10^{-6}\) on subset of \texttt{sslsq} problems with \( m, n \leq \num{20 000}\). All 32 problems were solved by at least one solver.}
    \label{fig:admm-time-relative-sslsq-eps1e-6}
\end{figure}

We use 32 Huber \cite{huberRobustEstimationLocation1992} and Lasso \cite{tibshiraniRegressionShrinkageSelection1996} fitting problems derived from the SuiteSparse collection \cite{davisUniversityFloridaSparse2011} as in the \texttt{OSQP} paper \cite{stellatoOSQPOperatorSplitting2020}, with \(m, n \leq \num{20000}\).
Time profiles for \(\epsilon = 10^{-3}\) and \(\epsilon = 10^{-6}\) are shown in \Cref{fig:admm-time-relative-sslsq-eps1e-3,fig:admm-time-relative-sslsq-eps1e-6} respectively.

Since these problems are relatively easy, the additional time expended by the Krylov methods operating on two working vectors and orthogonalising Arnoldi vectors does not pay off for \(\epsilon = 10^{-3}\).
But at \(\epsilon = 10^{-6}\), AA struggles to solve the most difficult 10\% of problems, as its conditioning difficulties begin to manifest when FOM convergence slows down near optimality.
Here, the Krylov method's work in avoiding conditioning issues \textit{does} pay off and it outperforms the AA incumbents; presumably having identified optimal active sets, requesting higher accuracy is not a big challenge.
The Krylov profiles for \(\epsilon = 10^{-3}\) and \(\epsilon = 10^{-6}\) are qualitatively similar.

As with the \texttt{mpc} set, iteration-count profiles show a consistent advantage for the Krylov variants.



\section{Conclusion}
\label{sec:conclusion}

We proposed a novel Krylov acceleration method for FOMs on QPs.
While equivalent in a sense to Anderson's, our method exploits the Arnoldi process at the heart of GMRES to avoid the conditioning issues which often plague the former in regions of slower convergence.

We included a brief discussion of a broadly applicable safeguarding rule with theory guarantees adapted from existing literature.
Further, we considered computational performance in sparse matrix-vector products and linear system solutions, the relevant primitives which our method has to perform for two working vectors simultaneously, in contrast with Anderson acceleration's requirements.

With thorough benchmarks on QPs arising in model predictive control and statistical learning, we demonstrated the performance of our method on an ADMM baseline compared to an efficient Anderson acceleration module used by the \texttt{COSMO} conic solver.

We found that our Krylov method is stronger across the board in terms of iterations.
In many cases it also performs better in wall time terms, particularly for optimal control problems and for problems solved to higher accuracy.

The general-operator formulation of \Cref{sec:fom-linearisation} extends naturally to FOMs for conic programs beyond QPs, where the FOM operator \(T\) involves projections onto more general cones (\eg second-order or positive semidefinite).
These projections are differentiable almost everywhere, so the Jacobian \(L = DT(u_k)\) and the linearised system \eqref{eq:linearised-fp-system} are well-defined at generic points.
Unlike in the QP case, where projections onto the nonnegative orthant are piecewise-affine and the linearisation is everywhere exact, general cone projections introduce genuinely approximate linearisations, so that the Arnoldi relation \eqref{eq:arnoldi-relation} only holds exactly under stronger conditions.
The safeguarding framework of \Cref{sec:krylov-adapted-safeguard} protects convergence regardless of proposal quality, making this extension a natural direction for future work.

\balance


\bibliographystyle{IEEEtran}
\bibliography{IEEEabrv,DPhil}


\end{document}